# A GENERALIZATION OF HALPERN´s ITERATION: PRELIMINARY RESULTS

M. De la Sen

**IIDP**. Faculty of Science and Technology.

**University of the Basque Country.** Campus of Leioa (Bizkaia). Aptdo. 644- Bilbao, SPAIN

email: manuel.delasen@ehu.es

**Abstract**. A generalization of Halpern´s iteration is investigated on a compact convex subset of a smooth Banach space. The modified iteration process consists of a combination of a viscosity term, an external sequence, a continuous nondecreasing function of a distance of points of an external sequence, which is not necessarily related to the solution of Halpern´s iteration, a contractive mapping, and a non-expansive one. The sum of the real coefficient sequences of four of the above terms is not required to be unity at each sample but it is assumed to converge asymptotically to unity. Halpern´s iteration solution is proven to converge strongly to a unique fixed point of the asymptotically nonexpansive mapping.

## 1. Introduction

Fixed-point theory is a powerful tool for investigating the convergence of the solutions of iterative discrete processes or that of the solutions of differential equations to fixed points in appropriate convex compact subsets of complete metric spaces or Banach spaces, in general, [1-12]. A key point is that the equations under study be driven by contractive maps or at least by asymptotically nonexpansive maps. By that reason, the fixed-point formalism is useful in stability theory to investigate the asymptotic convergence of the solution to stable attractors which being stable equilibrium points. The uniqueness of the fixed point is not required in the most general context although it can be sometimes suitable provided that only one such a point exists in some given problem. Therefore, the theory is useful for stability problems subject to multiple stable equilibrium points. Compared to Lyapunov´s stability theory, it may be a more powerful tool in cases when searching a Lyapunov functional is a difficult task or when there exist multiple equilibrium points, [1], [12]. Furthermore, it is not easy to obtain the value of the equilibrium points from that of the Lyapunov functional in the case that the last one is very involved. A generalization of the contraction principle in metric spaces by using continuous non-decreasing functions subject to an inequality-type constraint has been performed in [2]. The concept of n- times reasonable expansive mapping in a complete metric space is defined in [3] and proven to possess a fixed point. In [5], the T-stability of Picard´s iteration is investigated with T being a self-mapping of X where (X, d) is a complete metric space. The concept of T-stability is set as follows: if a solution sequence converges to an existing fixed point of T, then the error in terms of distance of any two consecutive values of any solution generated by Picard´s iteration converges asymptotically to zero. On the other hand, an important effort has been devoted to the investigation of the Halpern´s iteration scheme and many associate extensions during the last decades (see, for instance, [4], [6], [9] and [10]). The basic Halpern´s iteration is driven by an external sequence plus a contractive mapping whose two associate coefficient sequences sum unity for all samples, [9]. Recent extensions of the Halpern´s iteration to viscosity iterations have been proposed in [4] and [6]. In the first reference, a viscosity –type term is added as extra forcing term to the basic external sequence of Halpern´s scheme. In the second one, the external driving term is replaced with two ones, namely, a viscosity- type term plus an asymptotically nonexpansive mapping taking values on a left-reversible semigroup of asymptotically nonexpansive Lipschitzian mappings on a compact convex subset C of the Banach space X. The final iteration process investigated in [6] consists of three forcing terms, namely, a contraction on C, an asymptotically nonexpansive Lipschitzian mapping taking values in a left- reversible semigroup of mappings from a subset of that of bounded functions on its dual. It is proven that the solution converges to a unique common fixed point of all the set asymptotic nonexpansive mappings for any initial conditions on C. The objective of this paper is to investigate further generalizations for Halpern´s iteration process via fixed-point theory by using two more driving terms, namely, an external one taking values on C plus a nonlinear term given by a continuous nondecreasing function , subject to an inequality-type constraint as proposed in [2], whose argument is the distance between pairs of points of sequences in certain complete metric space which are not necessarily directly related to the sequence solution taking values in the subset C of the Banach space X. Another generalization point is that the sample-by-sample sum of the scalar coefficient sequences of all the driving terms is not necessarily unity but it converges asymptotically to unity.

## 2. Stability and boundedness properties of a viscosity-type difference equation

In this section a real difference equation scheme is investigated from a stability point of view by also discussing the existence of stable limiting finite points. The structure of such an iterative scheme supplies the structural basis for the general viscosity iterative scheme later discussed formally in Section 4 in the light of contractive and asymptotically nonexpansive mappings in compact convex subsets of Banach spaces. The following well–known iterative scheme is investigated for an iterative scheme which generates real sequences:

**Theorem 2.1.** Consider the difference equation:

$$x_{k+1} = \beta_k x_k + (1-\beta_k) z_k \tag{2.1}$$

such that the error sequence $\{e_k := x_k - z_k\}$ is generated by

$$e_{k+1} = \beta_k e_k - \tilde{z}_{k+1} \tag{2.2}$$

$\forall k \in \mathbf{Z}_{0+} := \mathbf{N} \cup \{0\}$, where $\tilde{z}_k := z_{k+1} - z_k$

Assume that $x_0$ and $z_0$ are bounded real constants and $0 \le \beta_k < 1$; $\forall k \in \mathbf{Z}_{0+}$. Then, the following properties hold:

**(i)** The real sequences $\{x_k\}$, $\{z_k\}$ and $\{e_k\}$ are uniformly bounded if $0 \le e_k \le \dfrac{2 x_k}{1-\beta_k}$ if $x_k > 0$ and $\dfrac{2 x_k}{1-\beta_k} \le e_k \le 0$ if $x_k \le 0$; $\forall k \in \mathbf{Z}_{0+}$. If, furthermore, $0 < e_k < \dfrac{2 x_k}{1-\beta_k}$ if $x_k > 0$ and $\dfrac{2 x_k}{1-\beta_k} < e_k \le 0$ if $x_k \le 0$, with $e_k = 0$ if and only if $x_k = 0$; $\forall k \in \mathbf{Z}_{0+}$, then the sequences $\{x_k\}$, $\{z_k\}$ and $\{e_k\}$ converge asymptotically to the zero equilibrium point as $k \to \infty$ and $\{|x_k|\}$ is monotonically decreasing.

**(ii)** Let the real sequence $\{\ell_k\}$ be defined by $\ell_k := \dfrac{\tilde{z}_{k+1}}{e_k} = \dfrac{z_{k+1} - z_k}{x_k - z_k}$ if $x_k \ne z_k$ and $\ell_k = 1$ if $x_k = z_k$ (what implies that $z_{k+1} = x_{k+1} = x_k = z_k$ from (2.1) and $\ell_k = 1$). Then, $\{e_k\}$ is uniformly bounded if $\ell_k \in [\beta_k - 1, 1 + \beta_k]$; $\forall k \in \mathbf{Z}_{0+}$. If, furthermore, $\ell_k \in (\beta_k - 1, 1 + \beta_k)$; $\forall k \in \mathbf{Z}_{0+}$ then $e_k \to 0$ as $k \to \infty$.

**(iii)** Let $x_0 \geq 0$ and $\{z_k\}$ a positive real sequence (i.e. all its elements are nonnegative real constants). Define $\ell_k := \frac{\tilde{z}_{k+1}}{e_k}$ if $x_k \neq z_k$ and $\ell_k = 1$ if $x_k = z_k$. Then, $\{x_k\}$ is a positive real sequence and $\{e_k\}$ is uniformly bounded if $\ell_k \in [0, 1 - \beta_k]$; $\forall k \in \mathbf{Z}_{0+}$. If, furthermore, $\ell_k \in (0, 1 - \beta_k)$; $\forall k \in \mathbf{Z}_{0+}$ then $e_k \to 0$ as $k \to \infty$.

**(iv)** If $|\beta_k| \leq 1$; $\forall k \in \mathbf{Z}_{0+}$ and $\sum_{k=0}^{\infty} |z_k| < \infty$ then $|x_k| < \infty$; $\forall k \in \mathbf{Z}_{0+}$.

If $|\beta_k| \leq \beta < 1$ and $|z_k| < \infty$; $\forall k \in \mathbf{Z}_{0+}$ then $|x_k| < \infty$; $\forall k \in \mathbf{Z}_{0+}$.

If $|\beta_k| \leq \beta < 1/(1 + 2\beta_0)$ and $|z_k| \leq \beta_0 |x_k| < \infty$; $\forall k \in \mathbf{Z}_{0+}$ for some $\beta_0 \in \mathbf{R}_+ := \{z \in \mathbf{R} : z > 0\}$, with $\mathbf{R}_{0+} := \{z \in R : z \geq 0\} = \mathbf{R}_+ \cup \{0\}$, then $|x_k| < \infty$; $\forall k \in \mathbf{Z}_{0+}$ and $x_k \to 0$ as $k \to \infty$.

**(v)** (Corollary to Venter's theorem, [7]). Assume that $\beta_k \in [0, 1]$; $\forall k \in \mathbf{Z}_{0+}$, $(1 - \beta_k) \to 0$ as $k \to \infty$ and $\sum_{j=0}^{k} (1 - \beta_j) \to \infty$ (what imply $\beta_k \to 1$ as $k \to \infty$ and the sequence $\{\beta_k\}$ has only a finite set of unity values). Assume also that $x_0 \geq 0$ and $\{z_k\}$ is a nonnegative real sequence with $\sum_{k=0}^{\infty} (1 - \beta_k) z_k < \infty$. Then $x_k \to 0$ as $k \to \infty$.

**(vi)** (Suzuki [8]- see also Saeidi [6]). Let $\{\beta_k\}$ be a sequence in $[0, 1]$ with $0 < \liminf_{k \to \infty} \beta_k \leq \limsup_{k \to \infty} \beta_k < 1$ and let $\{x_k\}$ and $\{z_k\}$ be bounded sequences. Then,
$$\limsup_{k \to \infty} (|z_{k+1} - z_k| - |x_{k+1} - x_k|) \leq 0.$$

**(vii)** (Halpern, [9]- see Hu, [4]). Let $z_k$ be $z_k = P x_k$; $\forall k \in \mathbf{Z}_{0+}$ in (2.1) subject to $x_0 \in C$, $\beta_k \in [0, 1]$; $\forall k \in \mathbf{Z}_{0+}$ with $P : C \to C$ being a nonexpansive self- mapping on C. Thus, $\{x_k\}$ converges weakly to a fixed point of P in the framework of Hilbert spaces endowed with the inner product $<x, Px>$; $\forall x \in X$, if $\beta_k = k^{-\beta}$ for any $\beta \in (0, 1)$.

**Proof**: **(i)** Direct calculations with (2.1) lead to:
$$x_{k+1}^2 - x_k^2 = (\beta_k^2 - 1) x_k^2 + (1 - \beta_k)^2 (x_k^2 + e_k^2 - 2 x_k e_k) + 2\beta_k (1 - \beta_k) x_k (x_k - e_k)$$
$$= (1 - \beta_k)^2 e_k^2 - 2(1 - \beta_k) x_k e_k = ((1 - \beta_k)^2 |e_k| - 2(1 - \beta_k) x_k \operatorname{sgn} e_k) |e_k| \text{ if } e_k \neq 0$$

so that $x_{k+1}^2 \leq x_k^2$ if $(1 - \beta_k)^2 e_k \operatorname{sgn} e_k \leq 2(1 - \beta_k) x_k \operatorname{sgn} e_k$, and equivalently, if $(1 - \beta_k) |e_k| \leq 2|x_k|$ and $e_k x_k = (x_k - z_k) x_k \geq 0$ with $e_k \neq 0$; and

$$x_{k+1}^2 - x_k^2 = 2\beta_k(1-\beta_k)x_k^2 \leq 0 \text{ if } e_k = x_k - z_k = 0$$

Thus, $x_{k+1}^2 \leq x_k^2 \leq x_0^2 < \infty$, $|e_k| \leq \dfrac{2|x_k|}{1-\beta_k} \leq \dfrac{2|x_0|}{1-\beta_k} < \infty$ and $|z_k| = \left|\dfrac{x_{k+1} - \beta_k x_k}{1-\beta_k}\right| \leq \dfrac{1+\beta_k}{1-\beta_k}|x_0| < \infty$;

$\forall k \in \mathbf{Z}_{0+}$. If, in addition, $(1-\beta_k)|e_k| < 2|x_k|$ and $e_k x_k = (x_k - z_k)x_k \geq 0$ with $e_k \neq 0$ then $x_k \to 0$ and is a monotonically decreasing sequence, $z_k \to 0$ and $e_k \to 0$ as $k \to \infty$. Property (i) has been proven

**(ii)** Direct calculations with (2.2) yield for $e_k \neq 0$

$$e_{k+1}^2 - e_k^2 = \left(\beta_k^2 - 1 + \ell_k^2 - 2\beta_k \ell_k\right) e_k^2 \leq 0 \text{ if } g(\ell_k) := \ell_k^2 - 2\beta_k \ell_k + \beta_k^2 - 1 \leq 0$$

Since $g(\ell_k)$ is a convex parabola $g(\ell_k) \leq 0$ for all $\ell \in \left[\ell_{k1}, \ell_{k2}\right]$ if real constants $\ell_{ki}$ exist such that $g(\ell_{ki}) = 0$; i=1,2. The parabola zeros are $\ell_{k1,2} = \beta_k \pm 1$ so that $e_{k+1}^2 \leq e_k^2 \leq e_0^2 < \infty$ if $\ell_k \in [\beta_k - 1, \beta_k + 1]$.

If $e_k = 0$ then $e_{k+1} = -\tilde{z}_{k+1} = z_k - z_{k+1} = x_{k+1} - z_{k+1} = e_k = 0$ with $\ell_k = 1$. Thus, $e_{k+1}^2 \leq e_k^2 \leq e_0^2 < \infty$ if $\ell_k \in [\beta_k - 1, \beta_k + 1]$; $\forall k \in \mathbf{Z}_{0+}$. If $\ell_k \in (\beta_k - 1, \beta_k + 1)$ then $e_k \to 0$ as $k \to \infty$. Property (ii) has been proven.

**(iii)** If $\{z_k\}$ is positive then $\{x_k\}$ is positive from direct calculations through (2.1). The second part follows directly from Property (ii) by restricting $\ell_k \in [0, \beta_k + 1]$ for uniform boundedness of $\{e_k\}$ and $\ell_k \in (0, \beta_k + 1)$ for its asymptotic convergence to zero in the case of nonzero $e_k$.

**(iv)** If $|\beta_k| \leq 1$; $\forall k \in \mathbf{Z}_{0+}$ and $\sum_{k=0}^{\infty} |z_k| < \infty$ then from recursive evaluation of (2.1):

$$|x_k| = \left|\prod_{j=0}^{k}[\beta_j]x_0 + \sum_{j=0}^{k}\prod_{\ell=j+1}^{k}[\beta_\ell](1-\beta_j)z_j\right| \leq |x_0| + \left|x_0 + \sum_{j=0}^{k} z_j\right| < \infty; \forall k \in \mathbf{Z}_{0+}$$

If, $|\beta_k| \leq \beta < 1$ and $|z_k| < \infty$; $\forall k \in \mathbf{Z}_{0+}$ then:

$$|x_k| \leq \left|\beta^k x_0\right| + \left|\sum_{j=0}^{k}\prod_{\ell=j+1}^{k}\beta^{k-\ell}(1-\beta_j)z_j\right| \leq \left|\beta^k x_0\right| + \dfrac{2}{1-\beta}\left(1-\beta^{k-1}\right)\max_{0 \leq j \leq k}|z_j| \leq |x_0| + \dfrac{2}{1-\beta}\max_{0 \leq j \leq k}|z_j| < \infty$$

; $\forall k \in \mathbf{Z}_{0+}$

If $|\beta_k| \leq \beta < 1/(1 + 2\beta_0)$ and $|z_k| \leq \beta_0 |x_k| < \infty$; $\forall k \in \mathbf{Z}_{0+}$ for some $\beta_0 \in \mathbf{R}_{0+} := \{0 \neq z \in \mathbf{R}_+\}$ then

$|x_{k+1}| \leq \beta |x_k| + 2\beta \beta_0 |x_k| \leq (1 + 2\beta_0)\beta |x_k| < |x_k|$ ; $\forall k \in \mathbf{Z}_{0+}$, thus $\{|x_k|\}$ is monotonically strictly decreasing so that it converges asymptotically to zero. □

Eq. (2.1) under the form

$$x_{k+1} = \beta_k x_k + (1 - \beta_k) P x_k \qquad (2.3)$$

with $x_0 \in C$ and $P: C \to C$ being a nonexpansive self-mapping on C under the weak or strong convergence conditions of Theorem 2.1 (vii)–(viii) is known as Halpern`s iteration, [4] which is a particular case of the generalized viscosity iterative scheme studied in the subsequent sections. Theorem 2.1 (vi) extends stability Venter´s theorem which is useful in recursive stochastic estimation theory when investigating the asymptotic expectation of the norm-squared parametrical estimation error, [7]. Note that the stability result of this section has been derived by using the discrete Lyapunov´s stability theorem with the Lyapunov´s sequence $\{V_k := x_k^2\}$ what guarantees global asymptotic stability to the zero equilibrium point if it is strictly monotonically decreasing on $\mathbf{R}_+$ and to global stability (stated essentially in terms of uniform boundedness of the sequence $\{x_k\}$) if it is monotonically decreasing on $\mathbf{R}_+$. The links between Lyapunov´s stability and fixed point theory are clear (see, for instance, [1-2]). However, fixed point theory is a more powerful tool uncertain problems since it copes more easily with the existence of multiple stable equilibrium points and with nonlinear mappings. Note that the results of Theorem 2.1 may be further formalized in the context of fixed-point theory by defining a complete metric space $(\mathbf{R}, d)$, respectively, $(\mathbf{R}_{0+}, d)$ for the particular results being applicable to a positive system under nonnegative initial conditions, with the Euclidean metrics defined by $d(x_k, z_k) = |x_k - z_k|$

## 3. Some definitions and background as preparatory tools for Section 4

The four subsequent definitions and then used in the results established and proven in Section 4.

**Definition 3.1**. *S* is a left reversible semigroup if $aS \cap bS \neq \emptyset$ ; $\forall a, b \in S$. □

It is possible to define a partial preordering relation "$\prec$" by $a \prec b \Leftrightarrow aS \supset bS$ ; $\forall a, b \in S$ for any semigroup *S*. Thus, $\exists c = aa' = bb' \in S$, for some existing a' and b' $\in S$, such that $aS \cap bS \supseteq cS \Rightarrow (a \prec c \land b \prec c)$ if *S* is left-reversible. The semigroup *S* is said to be left-amenable if it has a left-invariant mean and it is then left-reversible, [6], [13].

**Definition 3.2**. $\mathbf{S} := \{T(s) : s \in S\}$ is said to be a representation of a left reversible semigroup *S* as Lipschitzian mappings on C if $T(s)$ is a Lipschitzian mapping on C with Lipschitz constant $k(s)$ and, furthermore, $T(st) = T(s)T(t)$; $\forall s, t \in S$.

The representation $\mathbf{S} := \{T(s): s \in S\}$ may be non expansive, asymptotically nonexpansive, contractive and asymptotically contractive according to Definition 3.3 and Definitions 3.4 which follow:

**Definition 3.3.** A representation $\mathbf{S} := \{T(s): s \in S\}$ of a left reversible semigroup $S$ as Lipschitzian mappings on C, a nonempty weakly compact convex subset of X, with Lipschitz constants $\{k(s): s \in S\}$ is said to be a nonexpansive (respectively, asymptotically nonexpansive, [6]) semigroup on C if it holds the uniform Lipschitzian condition $k(s) \leq 1$ (respectively, $\lim_S k(s) \leq 1$) on the Lipschitz constants. □

**Definitions 3.4.** A representation $\mathbf{S} := \{T(s): s \in S\}$ of a left reversible semigroup $S$ as Lipschitzian mappings on C with Lipschitz constants $\{k(s): s \in S\}$ is said to be a contractive (respectively, asymptotically contractive) semigroup on C if it holds the uniform Lipschitzian condition $k(s) \leq \delta < 1$ (respectively, $\lim_S k(s) \leq \delta < 1$) on the Lipschitz constants. □

The iteration process (3.1) is subject to a forcing term generated by a set of Lipschitzian mappings $\mathbf{S} \ni T(\mu_k): Z^* \times C \to C$ where $\{\mu_k\}$ is a sequence of means on $Z \subset \ell^\infty(S)$, with the subset Z (defined in Definition 3.5 below) containing unity, where $\ell^\infty(S)$ is the Banach space of all bounded functions on $S$ endowed with the supremum norm, such that $\mu_k: Z \to Z^*$ where $Z^*$ is the dual of $Z$.

**Definition 3.5.** The real sequence $\{\mu_k\}$ is a sequence of means on Z if $\|\mu_k\| = \mu_k(1) = 1$. □

Some particular characterizations of sequences of means to be invoked later on in the results of Section 4 are now given in the definitions which follow:

**Definitions 3.6.** The sequence of means $\{\mu_k\}$ on $Z \subset \ell^\infty(S)$ is:

1. Left -invariant if $\mu(\ell_s f) = \mu(f)$; $\forall s \in S$, $\forall f \in Z$, $\forall \mu \in \{\mu_k\}$ in $Z^*$ for $\ell_s \in \ell^\infty(S)$.

2. Strongly left- regular if $\lim_\alpha \|\ell_s^* \mu_\alpha - \mu_\alpha\| = 0$, $\forall s \in S$, where $\ell_s^*$ is the adjoint operator of $\ell_s \in \ell^\infty(S)$ defined by $(\ell_s f)(t) = f(st)$; $\forall t \in S$, $\forall f \in \ell^\infty(S)$. □

Parallel definitions follow for right-invariant and strongly right-amenable sequences of means. Z is said to be left (respectively right)-amenable if it has a left (respectively right)–invariant mean. A general viscosity iteration process considered in [6] is the following:

$$x_{k+1} = \alpha_k f(x_k) + \beta_k x_k + \gamma_k T(\mu_k) x_k; \quad \forall k \in \mathbf{Z}_{0+} \tag{3.1}$$

where:

- the real sequences $\{\alpha_k\}$, $\{\beta_k\}$ and $\{\gamma_k\}$ have elements in $(0, 1)$ of sum being identity, $\forall k \in \mathbf{Z}_{0+}$

- $\mathbf{S} := \{T(s): s \in S\}$ is a representation of a left reversible semigroup with identity $S$ being asymptotically nonexpansive, on a compact convex subset C of a smooth Banach space, with respect to a left-regular sequence of means defined on an appropriate invariant subspace of $\ell^\infty(S)$

- f is a contraction on C

It has been proven that the solution of the sequence converges strongly to a unique common fixed point of the representation $\mathbf{S}$ which is the solution of a variational inequality, [6]. The viscosity iteration process (3.1) generalizes that proposed in [13] for $\alpha_k = 0$ and $\gamma_k = 1 - \beta_k$ and also that proposed in [14-15] with $\beta_k = 0$, $\gamma_k = 1 - \beta_k$ and $T(\mu_k) = T$; $\forall k \in \mathbf{Z}_{0+}$. Halpern´s iteration is obtained by replacing $\gamma_k T(\mu_k) \to (1 - \alpha_k)u$ and $\beta_k = 0$ in (3.1) by using the formalism of Hilbert spaces, $\forall k \in \mathbf{Z}_{0+}$ (see, for instance, [4], [9], [10]). It has been proven the weak convergence of the sequence $\{x_k\}$ to a fixed point of T for any given u, $x_0 \in \mathbf{C}$ if $\alpha_k = k^{-\alpha}$ for $\alpha \in (0,1)$, [9]; and also proven to converge strongly to one such a point if $\alpha_k \to 0$ and $\frac{\alpha_{k+1} - \alpha_k}{\alpha_{k+1}^2} \to 0$ as $k \to \infty$, and $\sum_{k=0}^{\infty} \alpha_k = +\infty$, [10]. On the other hand, note that if $\alpha_k = 0$ and $\gamma_k = 1 - \beta_k$ and $z_k = T(\mu_k)x_k$ with $x_k \in \mathbf{R}$; $\forall k \in \mathbf{Z}_{0+}$ then the resulting particular iteration process (3.1) becomes the difference equation (2.1) discussed in Theorem 2.1 from a stability point of view provided that the boundedness of the solution is ensured on some convex compact set $\mathbf{C} \subset \mathbf{R}$; $\forall k \in \mathbf{Z}_{0+}$.

**4. Boundedness and convergence properties of a more general difference equation**

The viscosity iteration process (3.1) is generalized in this section by including two more forcing terms not being directly related to the solution sequence. One of them being dependent on a nondecreasing distance-valued function related to a complete metric space while the other forcing term is governed by an external sequence $\{\delta_k r\}$. Furthermore the sum of the four terms of the scalar sequences $\{\alpha_k\}$, $\{\beta_k\}$ and $\{\gamma_k\}$ and $\{\delta_k\}$ at each sample is not necessarily unity but it is asymptotically convergent to unity.

The following generalized viscosity iterative scheme, which is a more general difference equation than (4.1), is considered in the sequel:

$$x_{k+1} = \alpha_k f(x_k) + \beta_k x_k + \gamma_k T(\mu_k) x_k + \left( \sum_{i=1}^{s_k} v_{ik} \varphi_i \left( d(\omega_k, \omega_{k-p}) \right) + \delta_k r \right); \; \forall k \in \mathbf{Z}_{0+} \quad (4.1)$$

$\forall x_0 \in \mathbf{C}$ for a sequence of given finite numbers $\{s_k\}$ with $s_k \in \mathbf{Z}_{0+}$ ( if $s_k = 0$ then the corresponding sum is dropped off) which can be rewritten as (2.1) if $0 < \beta_k < 1$; $\forall k \in \mathbf{Z}_{0+}$ (except possibly for a finite number of values of the sequence $\{\beta_k\}$ what implies $0 < \liminf_{k \to \infty} \beta_k \leq \limsup_{k \to \infty} \beta_k < 1$ ) by defining the sequence

$$z_k = \frac{1}{1 - \beta_k} \left( \alpha_k f(x_k) + \gamma_k T(\mu_k) x_k + \left( \sum_{i=1}^{s_k} v_{ik} \varphi_i \left( d(\omega_k, \omega_{k-p}) \right) + \delta_k r \right) \right) \quad (4.2)$$

with $x_0 \in \mathbf{C}$, where:

- $\{\mu_k\}$ is a strongly left regular sequence of means on $Z \subset \ell^\infty(S)$, i.e. $\mu_k \in Z^*$. See Definition 3.5 and Definition 3.5.

- $S$ is a left reversible semigroup represented as Lipschitzian mappings on C by $\mathbf{S} := \{T(s) : s \in S\}$.

The iterative scheme is subject to the following assumptions:

**Assumptions 4.1**

1. $\{\alpha_k\}$, $\{\gamma_k\}$ and $\{\delta_k\}$ are real sequences in $[0,1]$, $\{\beta_k\}$ is a real sequence in $[0,1)$, $\{\nu_{ik}\}$ are sequences in $\mathbf{R}_{0+}$; $\forall i \in \overline{k} := \{1, 2, ..., k\}$ for some given $k \in \mathbf{Z}_+ \equiv \mathbf{N} := \mathbf{Z}_{0+} \setminus \{0\}$ and $r \in \mathbf{R}$

2. $\lim_{k \to \infty} \alpha_k = \lim_{n \to \infty} \delta_k = 0$, $\liminf_{k \to \infty} \gamma_k > 0$

3. $\lim_{k \to \infty} \sum_{j=1}^{k} \alpha_j = \infty$, $\lim_{k \to \infty} \sum_{j=1}^{k} \delta_j < \infty$

4. $0 < \liminf_{k \to \infty} \beta_k \le \limsup_{k \to \infty} \beta_k < 1$

5. $\alpha_k + \beta_k + \gamma_k + \delta_k = 1 + (1 - \beta_k)\varepsilon_k$; $\forall k \in \mathbf{Z}_{0+}$ with $\{\varepsilon_k\}$ being a bounded real sequence satisfying $\varepsilon_k \ge \dfrac{1}{\beta_k - 1}$ and $\lim_{k \to \infty} \varepsilon_k = 0$

6. $f$ is a contraction on a nonempty compact convex subset C, of diameter $d_C = \operatorname{diam} C := \sup \{\|x - y\| : x, y \in C\}$, of a Banach space X, of topological dual $X^*$, which is smooth, i.e. its normalized duality mapping $J: X \to 2^{X^*} \subset X^*$ from X into the family of nonempty (by the Hahn-Banach theorem [6], [11]), weak-star compact convex subsets of $X^*$, defined by

$$J(x) := \left\{x^* \in X^* : x^*(x) = <x, x^*> = \|x^*\|^2 = \|x\|^2\right\} \subset X^*, \ \forall x \in X$$

is single-valued

7. The representation $\mathbf{S} := \{T(s) : s \in S\}$ of the left-reversible semigroup $S$ with identity, is asymptotically nonexpansive on C (see Definition 3.3) with respect to $\{\mu_k\}$, with $\mu_k \in Z^*$ which is strongly left regular so that it fulfils $\lim_{k \to \infty} \|\mu_{k+1} - \mu_k\| = 0$

8. $\limsup_{k \to \infty} \sup_{x, y \in C} \left(\|T(\mu_k)x - T(\mu_k)y\| - \|x - y\|\right) / \min(\alpha_k, \delta_k) \le 0$

9. $(W, d)$ is a complete metric space and $Q: W \to W$ is a self-mapping satisfying the inequality

$\varphi_i(d(Qy, Qz)) \le \varphi_i(d(y, z)) - \phi_i(d(y, z))$; $\forall y, z \in W$

where $\varphi_i, \phi_i \in \mathbf{R}_{0+} \to \mathbf{R}_{0+}$; $\forall i \in \overline{k}$ are continuous monotone nondecreasing functions satisfying $\varphi_i(t) = \phi_i(t) = 0$ if and only if $t = 0$; $\forall i \in \overline{k}$

10. $\{\omega_k\}$ is a sequence in W generated as $\omega_{k+1} = Q\omega_k$, $k \in \mathbf{Z}_{0+}$ with $\omega_0 \in W$ and $p \in \mathbf{Z}_+$ is a finite given number □

Note that Assumption 4.1.4 is stronger than the conditions imposed on the sequence $\{\beta_k\}$ in Theorem 2.1 for (2.1). However, the whole viscosity iteration is much more general than the iterative equation (2.1). Three generalizations compared to existing schemes of this class are that an extra coefficient sequence $\{\delta_k\}$ is added to the set of usual coefficient sequences and that the exact constraint for the sum of coefficients $\alpha_k + \beta_k + \gamma_k + \delta_k$ being unity for all k is replaced by a limit-type constraint $\alpha_k + \beta_k + \gamma_k + \delta_k \to 1$ as $k \to \infty$ while during the transient such a constraint can exceed unity or be below unity at each sample (see Assumption 4.1.5). Another generalization is the inclusion of a nonnegative term with generalized contractive mapping $Q: W \to W$ involving another iterative scheme evolving on another, and in general distinct, complete metric space $(W, d)$ (see Assumptions 4.1.9 and 4.1.10). Some boundedness and convergence properties of the iterative process (4.1) are formulated and proven in the subsequent result.

**Theorem 4.2**. The difference iterative scheme (4.1), and equivalently the difference equation (2.1) subject to (4.2), possesses the following properties under Assumptions 4.1:

**(i)** $\max\left(\sup_{k \in \mathbf{Z}_{0+}} |x_k|, \sup_{k \in \mathbf{Z}_{0+}} |T(\mu_k) x_k|\right) < \infty$ ; $\forall x_0 \in C$. Also, $\|x_k\| < \infty$ and $\|T(\mu_k) x_k\| < \infty$ for any norm defined on the smooth Banach space X and there exists a nonempty bounded compact convex set $C_0 \subseteq C \subset X$ such that the solution of (4.2) is permanent in $C_0$, $\forall k \geq k_0$ and some sufficiently large finite $k_0 \in \mathbf{Z}_{0+}$ with $\max_{k \geq k_0}\left(\|x_k\|, \|T(\mu_k) x_k\|\right) \leq d_{C_0} := \text{diam } C_0$.

**(ii)** $\lim_{k \to \infty} \|T(\mu_k) x_k - x_k\| = 0$ and $x_k \to z_k \to \dfrac{\gamma_k T(\mu_k) x_k}{1 - \beta_k} \to T(\mu_k) x_k \to x^* \in C_0$ as $k \to \infty$.

**(iii)** $\infty > \left|x^* - x_0\right| = \left|\lim_{k \to \infty} \sum_{j=0}^{k}(x_{j+1} - x_j)\right|$

$= \left|\sum_{j=0}^{\infty}\left(\alpha_j f(x_j) + (\beta_j - 1)x_j + \gamma_j T(\mu_j) x_j + \left(\sum_{i=1}^{s_j} v_{ij} \varphi_i(d(\omega_j, \omega_{j-p})) + \delta_j r\right)\right)\right|$

**(iv)** Assume that $\{x_k\} \in C$ such that each sequence element $x_k \in \mathbf{R}_{0+}^m$ ( the first closed orthant of $\mathbf{R}^m$) ; $\forall k \in \mathbf{Z}_{0+}$ , for some $m \in \mathbf{Z}_+$ so that (4.1) is a positive viscosity iteration scheme. Then,

   **(iv.1)** $\{x_k\}$ is a nonnegative sequence (i.e. all its components are nonnegative $\forall k \geq 0, \forall x_0 \in C$), denoted as $x_k \geq 0$; $\forall k \geq 0$.

   **(iv.2)** Property (i) holds for $C_0 \subseteq C$ and Property (ii) also holds for a limiting point $x^* \in C_0$.

   **(iv.3)** Property (iii) becomes:

$\infty > \left|x^* - x_0\right|$

$= \left|\sum_{j=0}^{\infty}\left(\alpha_j f(x_j) + \gamma_j T(\mu_j) x_j + \left(\sum_{i=1}^{s_j} v_{ij} \varphi_i(d(\omega_j, \omega_{j-p})) + \delta_j r\right)\right) - \sum_{j=0}^{\infty}\left((1 - \beta_j) x_j\right)\right|$

what implies that either

$$\sum_{j=0}^{\infty}\left(\alpha_j f(x_j) + \gamma_j T(\mu_j) x_j + \left(\sum_{i=1}^{s_j} \nu_{ij}\varphi_i\left(d\left(\omega_j, \omega_{j-p}\right)\right) + \delta_j r\right)\right) < \infty \text{ and } \sum_{j=0}^{\infty}\left((1-\beta_j) x_j\right) < \infty$$

or

$$\limsup_{k\to\infty} \sum_{j=0}^{k}\left(\alpha_j f(x_j) + \gamma_j T(\mu_j) x_j + \left(\sum_{i=1}^{s_j} \nu_{ij}\varphi_i\left(d\left(\omega_j, \omega_{j-p}\right)\right) + \delta_j r\right)\right) = \infty, \text{ and}$$

$$\limsup_{k\to\infty} \sum_{j=0}^{\infty}\left((1-\beta_j) x_j\right) = \infty$$

**Proof**: From (2.5) and substituting the real sequence $\{\gamma_k\}$ from the constraint Assumption 4.1.5:

$$z_{k+1} - z_k = \frac{1}{1-\beta_{k+1}}\left(\alpha_{k+1} f(x_{k+1}) + \gamma_{k+1} T(\mu_{k+1}) x_{k+1} + \left(\sum_{i=1}^{s_{k+1}} \nu_{i,k+1}\varphi_i\left(d\left(\omega_{k+1}, \omega_{k+1-p}\right)\right) + \delta_{k+1} r\right)\right)$$

$$- \frac{1}{1-\beta_k}\left(\alpha_k f(x_k) + \gamma_k T(\mu_k) x_k + \left(\sum_{i=1}^{s_k} \nu_{i,k}\varphi_i\left(d\left(\omega_k, \omega_{k-p}\right)\right) + \delta_k r\right)\right)$$

$$= \frac{1}{1-\beta_{k+1}}\left(\alpha_{k+1} f(x_{k+1}) + \left(1 + (1-\beta_{k+1})\varepsilon_{k+1} - \alpha_{k+1} - \beta_{k+1} - \delta_{k+1}\right) T(\mu_{k+1}) x_{k+1} + \left(\sum_{i=1}^{s_{k+1}} \nu_{i,k+1}\varphi_i\left(d\left(\omega_{k+1}, \omega_{k+1-p}\right)\right) + \delta_{k+1} r\right)\right)$$

$$- \frac{1}{1-\beta_k}\left(\alpha_k f(x_k) + \left(1 + (1-\beta_k)\varepsilon_k - \alpha_k - \beta_k - \delta_k\right) T(\mu_k) x_k + \left(\sum_{i=1}^{s_k} \nu_{i,k}\varphi_i\left(d\left(\omega_k, \omega_{k-p}\right)\right) + \delta_k r\right)\right)$$

$$= \left(1 - \frac{\alpha_{k+1} + \delta_{k+1}}{1-\beta_{k+1}} + \varepsilon_{k+1}\right) T(\mu_{k+1}) x_{k+1} - \left(1 - \frac{\alpha_k + \delta_k}{1-\beta_k} + \varepsilon_k\right) T(\mu_k) x_k$$

$$+ \frac{\alpha_{k+1}}{1-\beta_{k+1}} f(x_{k+1}) - \frac{\alpha_k}{1-\beta_k} f(x_k) + \left(\frac{\delta_{k+1}}{1-\beta_{k+1}} - \frac{\delta_k}{1-\beta_k}\right) r$$

$$+ \frac{1}{1-\beta_{k+1}}\left(\sum_{i=1}^{s_{k+1}} \nu_{i,k+1}\varphi_i\left(d\left(\omega_{k+1}, \omega_{k+1-p}\right)\right)\right) - \frac{1}{1-\beta_k}\left(\sum_{i=1}^{s_k} \nu_{i,k}\varphi_i\left(d\left(\omega_k, \omega_{k-p}\right)\right)\right) \quad (4.3)$$

Thus,

$$\|z_{k+1} - z_k\| \leq \|T(\mu_{k+1}) x_{k+1} - T(\mu_k) x_k\| + \left\|\left(\frac{\alpha_{k+1} + \delta_{k+1}}{1-\beta_{k+1}} + \varepsilon_{k+1}\right) T(\mu_{k+1}) x_{k+1} - \left(\frac{\alpha_k + \delta_k}{1-\beta_k} + \varepsilon_k\right) T(\mu_k) x_k\right\|$$

$$+ K_1\left(\alpha_k + \alpha_{k+1} + (\delta_k + \delta_{k+1})\right)|r| + K_2 \bar{s}\bar{\nu}\right) \; ; \; \forall k \geq k_0$$

$$\leq \|T(\mu_{k+1}) x_{k+1} - T(\mu_k) x_{k+1}\| + \|T(\mu_k) x_{k+1} - T(\mu_k) x_k\|$$

$$+ \left\|\left((\alpha_{k+1} + \delta_{k+1}) K_1 + \varepsilon_{k+1}\right) T(\mu_{k+1}) x_{k+1} - \left((\alpha_k + \delta_k) K_1 + \varepsilon_k\right) T(\mu_k) x_k\right\|$$

$$+ K\left((\alpha_k + \alpha_{k+1}) K_1 + (\delta_k + \delta_{k+1})\right)|r| + K_2 \bar{s}\bar{\nu}\right) \; ; \; \forall k \geq k_0$$

$$\leq \|T(\mu_{k+1})x_{k+1} - T(\mu_k)x_{k+1}\| + \|T(\mu_k)x_{k+1} - T(\mu_k)x_k\|$$

$$+ ((\alpha_k + \delta_k)K_1 + \varepsilon_k)\|(1+\rho_k)T(\mu_{k+1})x_{k+1} - T(\mu_k)x_k\|$$

$$+ K((\alpha_k + \alpha_{k+1})K_1 + (\delta_k + \delta_{k+1})|r| + K_2 \bar{s}\bar{v}) \ ; \ \forall k \geq k_0$$

$$\leq (1 + ((\alpha_k + \delta_k)K_1 + \varepsilon_k))(\|T(\mu_{k+1})x_{k+1} - T(\mu_k)x_{k+1}\|)$$

$$+ (1 + ((\alpha_k + \delta_k)K_1 + \varepsilon_k))\|T(\mu_k)x_{k+1} - T(\mu_k)x_k\|$$

$$+ ((\alpha_k + \delta_k)K_1 + \varepsilon_k)\|\rho_k T(\mu_{k+1})x_{k+1} - T(\mu_k)x_k\|$$

$$+ K((\alpha_k + \alpha_{k+1})K_1 + (\delta_k + \delta_{k+1})|r| + K_2 \bar{s}\bar{v}) \ ; \ \forall k \geq k_0 \qquad (4.4)$$

where $k_0 \in \mathbf{Z}_{0+}$ is an arbitrary finite sufficiently large integer, and

$$\bar{s} = \bar{s}(k_0) := \max_{k \geq k_0} s_k, \quad \bar{v} = \bar{v}(k_0) := \max_{k \geq k_0} \max_{i \in \bar{s}_k} v_{ik}$$

$$\rho_k := (\alpha_{k+1} + \delta_{k+1} - \alpha_k - \delta_k)K_1 + \varepsilon_{k+1} - \varepsilon_k \ ; \ \forall k \in \mathbf{Z}_{0+}$$

$$K := \frac{1}{1 - \limsup_{k \to \infty} \beta_k - \varepsilon_\beta} < \infty, \quad K_1 = K_1(x_0, k_0) := \sup_{k \geq k_0} |f(x_k)| \leq \sup_{x \in C} |f(x)| < \infty$$

$$\infty > K_2 = K_2(\omega_0, k_0) := 2\bar{s}(k_0)\bar{v}(k_0) \sup_{k \geq k_0} \max_{i \in \bar{s}_k} \varphi_i(d(\omega_k, \omega_{k-p})) \to 0 \text{ as } k_0 \to \infty$$

since the functions $\varphi_i$ are continuous on $\mathbf{R}_{0+}$ with $\varphi_i(0) = 0$ and $d(\omega_k, \omega_{k-p}) \to 0$ as $k \to \infty$, [2] with $\varepsilon_\beta > 0$ being prefixed and arbitrarily small. The constants $K, K_1$ and $K_2$ are finite for sufficiently large $k \in \mathbf{Z}_{0+}$ since $\limsup_{k \to \infty} \beta_k < 1$ (Assumption 4.1.4), f is a contraction on C (Assumption 4.1.6), and Q is a self-mapping on W satisfying Assumption 4.1.9. Since $\alpha_k \to 0$, $\delta_k \to 0$ and $\varepsilon_k \to 0$ as $k \to \infty$ from Assumptions 4.1.1 and 4.1.5 and $K_1$ is finite, $\rho_k \to 0$ as $k \to \infty$ and $|\rho_k| \leq \bar{\rho}(k_0)$; $\forall k \geq k_0$ being arbitrarily small since $k_0$ is arbitrarily large. Since from Assumption 4.1.7, **S** is an asymptotically nonexpansive semigroup on C, and $\alpha_k \to 0$, $\delta_k \to 0$ and $\varepsilon_k \to 0$ as $k \to \infty$:

$$(1 + ((\alpha_k + \delta_k)K_1 + \varepsilon_k))\|T(\mu_k)x_{k+1} - T(\mu_k)x_k\|$$

$$+ \left((\alpha_k + \delta_k)K_1 + \varepsilon_k\right)\left\|\rho_k T(\mu_{k+1})x_{k+1} - T(\mu_k)x_k\right\|$$

$$\leq (1+\varsigma_k)\|x_{k+1} - x_k\| + \xi_k, \quad \forall k \geq k_0 \tag{4.5}$$

with $\mathbf{R}_{0+} \ni \varsigma_k, \xi_k \to 0$ as $k \to \infty$. One gets from (4.5) into (4.4),

$$\|z_{k+1} - z_k\| \leq \left(1 + \left((\alpha_k + \delta_k)K_1 + \varepsilon_k\right)\right)\left(\|T(\mu_{k+1})x_{k+1} - T(\mu_k)x_{k+1}\|\right)$$

$$+ (1+\varsigma_k)\|x_{k+1} - x_k\| + \xi_k + K\left((\alpha_k + \alpha_{k+1})K_1 + (\delta_k + \delta_{k+1})|r| + \bar{s}\bar{v}K_2(\omega_0, k)\right); \quad \forall k \geq k_0$$

$$\tag{4.6}$$

what implies that

$$\limsup_{k \to \infty} \left(\|z_{k+1} - z_k\| - \|x_{k+1} - x_k\|\right) \leq \limsup_{k \to \infty} \left(\|z_{k+1} - z_k\| - \varsigma_k\|x_{k+1} - x_k\|\right)$$

$$\leq \limsup_{k \to \infty} \left(\left(1 + \left((\alpha_k + \delta_k)K_1 + \varepsilon_k\right)\right)\left(\|T(\mu_{k+1})x_{k+1} - T(\mu_k)x_{k+1}\|\right)\right.$$

$$\left. + \xi_k + K\left((\alpha_k + \alpha_{k+1})K_1 + (\delta_k + \delta_{k+1})|r| + \bar{s}\bar{v}K_2(\omega_0, k)\right)\right) = 0 \Rightarrow \lim_{k \to \infty}\|x_k - z_k\| = 0$$

$$\tag{4.7}$$

(see [8]) since $\|T(\mu_{k+1})x_{k+1} - T(\mu_k)x_{k+1}\| \to 0$ as $k \to \infty$ since $\{x_k\}$ is in C and $\{\mu_k\}$ is a strongly left regular sequence of means on X such that $\lim_{k \to \infty}\|\mu_{k+1} - \mu_k\| = 0$; and, furthermore, $\alpha_k \to 0$, $\delta_k \to 0$, $\varepsilon_k \to 0$, $\varsigma_k \to 0$, $\xi_k \to 0$ as $k \to \infty$ and $K_2(\omega_0, k) \to 0$ as $k \to \infty$. Thus, from (4.7) and using the above technical result in [8] for difference equations of the class (2.1) (see also [2]), it follows that:

$$\lim_{k \to \infty}\|x_{k+1} - x_k\| = \lim_{k \to \infty}(1-\beta_k)\|x_k - z_k\| = 0 \Rightarrow \lim_{k \to \infty}\|x_{k+1} - x_k\| = \lim_{k \to \infty}\|x_k - z_k\| = 0$$

$$\Rightarrow x_{k+1} \to x_k \to z_k \to \frac{\gamma_k T(\mu_k)x_k}{1-\beta_k} \text{ as } k \to \infty \tag{4.8}$$

since $0 < \liminf_{k \to \infty} \beta_k \leq \limsup_{k \to \infty} \beta_k < 1$ from (2.5) since $\alpha_k \to 0$, $\delta_k \to 0$ and $\varepsilon_k \to 0$ as $k \to \infty$. From (4.1):

$$x_{k+1} - x_k = \alpha_k f(x_k) + (1-\beta_k)(T(\mu_k)x_k - x_k) + ((1-\beta_k)\varepsilon_k - \alpha_k - \delta_k)T(\mu_k)x_k$$

$$+\left(\sum_{i=1}^{s_k} \nu_{ik}\varphi_i\left(d\left(\omega_k, \omega_{k-p}\right)\right)+\delta_k r\right); \quad \forall k \in \mathbf{Z}_{0+} \tag{4.9}$$

so that

$$\|T(\mu_k)x_k - x_k\| = \frac{1}{1-\beta_k}\left(\|x_{k+1} - x_k\| + \alpha_k\|f(x_k) - T(\mu_k)x_k\| + \left((1-\beta_k)\varepsilon_k - \alpha_k - \delta_k\right)\|T(\mu_k)x_k\|\right)$$

$$+\left(\sum_{i=1}^{s_k} \nu_{ik}\varphi_i\left(d\left(\omega_k, \omega_{k-p}\right)\right)+\delta_k|r|\right); \quad \forall k \in \mathbf{Z}_{0+} \tag{4.10}$$

Using Assumptions 4.1 and using (4.8) into (4.10) yields:

$$\lim_{k\to\infty} \|T(\mu_k)x_k - x_k\| = 0 \implies x_k \to T(\mu_k)x_k \text{ as } k \to \infty \tag{4.11}$$

since $\varphi_i\left(d\left(\omega_k, \omega_{k-p}\right)\right) \to 0$, $\alpha_k \to 0$, $\delta_k \to 0$, $\varepsilon_k \to 0$ as $k \to \infty$. Also, it follows that $x_k \to z_k \to \frac{\gamma_k T(\mu_k)x_k}{1-\beta_k} \to T(\mu_k)x_k$ as $k \to \infty$ from (4.8) and (4.11). Note that it has not been yet proven that the sequences $\{x_k\}$ and $\{z_k\}$ converge to a finite limit as $k \to \infty$ since it has not been proven that they are bounded. Thus, the four sequences $\{x_k\}, \{z_k\}, \left\{\frac{\gamma_k T(\mu_k)x_k}{1-\beta_k}\right\}$ and $\{T(\mu_k)x_k\}$ converge asymptotically to the same finite or infinite real limit. Proceed recursively with the solution of (4.1). Thus, for a given sufficiently large finite $n \in \mathbf{Z}_{0+}$ and $\forall k \in \mathbf{Z}_+$, one gets:

$$|x_{k+n}| = \left|\left(\prod_{i=n}^{k+n-1}[\beta_i]\right)x_n + \sum_{\ell=n}^{k+n-1}\left(\left(\prod_{j=\ell+1}^{k+n-1}[\beta_j]\right)\left(\alpha_\ell f(x_\ell) + \left((1-\beta_\ell)(1+\varepsilon_\ell)\varepsilon_\ell - \alpha_\ell - \delta_\ell\right)T(\mu_\ell)x_\ell\right)\right.\right.$$

$$\left.\left.+ \sum_{j=1}^{s_\ell}\nu_{j\ell}\varphi_j\left(d\left(\omega_\ell, \omega_{\ell-p}\right)\right)+\delta_\ell r\right)\right| \tag{4.12}$$

$$\leq \sigma^k M_n + \frac{1-\sigma^k}{1-\sigma}\left(\sup_{0\leq j\leq k+n-1}\rho_{j+n} + \left(\sup_{0\leq j\leq k+n-1}\lambda_{j+n}\right)\left(\sup_{0\leq j\leq k+n-1}\overline{M}_{j+n}\right)\right) \leq M_{k+n} \tag{4.13}$$

$\forall x_0 \in C_0$, for some positive real sequences $\{M_{j+n}\}$, $\{\rho_{j+n}\}$ and $\{\overline{M}_{j+n}\}$ satisfying $M_{j+n} \geq \sup_{0\leq i\leq j}|x_{i+n}|$ and $\overline{M}_{j+n} \geq \sup_{0\leq i\leq j}|T(\mu_{i+n})x_{i+n}|$, $\infty > \rho \geq \rho_{j+n} \to 0$ and $\infty > \lambda \geq \lambda_{j+n} \to 0$ as $j \to \infty$ with $\rho = \rho(n) > 0$ and $\lambda = \lambda(n) > 0$ being arbitrarily small for sufficiently large $n \in \mathbf{Z}_{0+}$, and

$$0 < \overline{\sigma} = \overline{\sigma}(n, n+1, \ldots, n+k-1) := 1 - \max_{n\leq j\leq n+k-1}\beta_j \leq \sigma := 1 - \limsup_{k\to\infty}\beta_k - b < 1$$

for sufficiently large $n \in \mathbf{Z}_{0+}$ and a sufficiently small $\mathbf{R}_+ \ni b = b(n) < 1 - \limsup_{k\to\infty}\beta_k \in (0,1)$ which exists from Assumptions 4.1.1 and 4.1.4. Note that the sequences $\{M_{j+n}\}$ and $\{\overline{M}_{j+n}\}$ may be chosen to satisfy $M_n \leq M_{j+n}$ and $\overline{M}_n \leq \overline{M}_{j+n}$; $\forall j \in \mathbf{Z}_{0+}$. Now, proceed by complete induction by assuming that $0 < \sup_{-n \leq j \leq k-1} \max(M_{j+n}, \overline{M}_{j+n}) \leq M < \infty$ for given sufficiently large $n \in \mathbf{Z}_{0+}$ and finite $k \in \mathbf{Z}_+$. Then, one gets from (4.13) that $0 < \sup_{-n \leq j \leq k} \max(M_{j+n}, \overline{M}_{j+n}, M_0) \leq M < \infty$ for any prescribed $M_0 \in \mathbf{R}_+$ if

$$\sigma^k M + \frac{1-\sigma^k}{1-\sigma}(\rho + \lambda M) \leq M \Leftrightarrow \frac{\rho}{M} + \lambda \leq 1 - \sigma \Leftrightarrow 0 < \sigma \leq 1 - \lambda - \frac{\rho}{M} \qquad (4.14)$$

with $\lambda = \lambda(n)$ and $\rho = \rho(n)$ which always holds for sufficiently large finite $n \in \mathbf{Z}_{0+}$ since $0 \leftarrow \max(\rho(n), \lambda(n)) \leq \frac{(1-\sigma)M}{M+1} < 1 - \sigma$ as $n \to \infty$. It has been proven by complete induction that the first part of Property (i) holds with the set $C_0$ being built such that $M = d_{C_0} = \text{diam } C_0$ for the given initial condition $x_0$. For a set of initial conditions $x_0 \in C_{00} \subset X$ with any set $C_{0in} \subset X$ convex and bounded, a common set $C_0$ might be defined for any initial condition of (4.1) in $C_{00}$ with a redefinition of the constant M as $M = \sup(M_{xo} : x_0 \in C_{00}) = d_{C_0} = \text{diam } C_0$. The second part of Property (i) follows for any norm on E from the property of equivalence of norms. Furthermore, the real sequences $\{x_k\}, \{z_k\}, \left\{\frac{\gamma_k T(\mu_k) x_k}{1-\beta_k}\right\}$ and $\{T(\mu_k) x_k\}$ converge strongly to a finite limit in $C_0$ since they are uniformly bounded so that Property (ii) has also been proven. Property (iii) follows directly from (4.1) and Property (ii). Property (iv.1) follows since $\{x_k\}$ is a nonnegative m-vector sequence provided that $x_0 \in \mathbf{R}_{0+}^m$ if $r \in \mathbf{R}_{0+}$ what follows from simple inspection of (4.1). Properties (iv.2)-(iv.3) follow directly from separating nonnegative positive and nonpositive terms in the right-hand side of the expression in Property (iii). □

The convergence properties of Theorem 4.2 (ii) are now related to the limits being fixed points of the asymptotically nonexpansive semigroup $\mathbf{S} := \{T(s) : s \in S\}$ which is the representation as Lipschitzian mappings on C of a left reversible semigroup $S$ with identity.

**Theorem 4.3**. The following properties hold:

**(i)** Let $F(\mathbf{S}) \in C$ be the set of fixed points of the asymptotically nonexpansive semigroup $\mathbf{S}$ on C. Then, the common strong limit $x^* \in C_0 \subseteq C$ of the sequences $\{x_k\}, \{z_k\}, \left\{\frac{\gamma_k T(\mu_k) x_k}{1-\beta_k}\right\}$ and $\{T(\mu_k) x_k\}$ in Theorem 4.2 (ii) is a fixed point of C located in $C_0$ and, thus, a stable equilibrium point of the iterative scheme (4.1) provided that $\text{diam } C_0$, and then $\text{diam } C$, is sufficiently large.

**(ii)** $F(\mathbf{S}) \subseteq C_0 \subseteq C$.

**Proof**: (i) Proceed by contradiction by assuming that $C_0 \ni x^* \notin F(\mathbf{S})$ so that there exists $\varepsilon_T \in \mathbf{R}_+$ such that

$$0 < \varepsilon_T \le \liminf_{k\to\infty} \|T(\mu_k)x^* - x^*\| \le \limsup_{k\to\infty} \|T(\mu_k)x^* - T(\mu_k)x_k\|$$

$$+ \left(\limsup_{k\to\infty} \|T(\mu_k)x_k - x_k\| + \limsup_{k\to\infty} \|x_k - x^*\|\right) = \limsup_{k\to\infty} \|T(\mu_k)x^* - T(\mu_k)x_k\|$$

$$\le \limsup_{k\to\infty} \|x_k - x^*\| = 0$$

since $\lim_{k\to\infty} \|T(\mu_k)x_k - x_k\| + \lim_{k\to\infty} \|x_k - x^*\| = 0$ where the above two limits exist and are zero from Theorem 3.2 (ii). Then, $x^* \in F(\mathbf{S})$, with $F(\mathbf{S})$ being nonempty since, at least one such finite fixed point exists in $C_0 \subseteq C$.

Property (ii) follows directly from Theorem 4.2 (iii)–(iv). □

**Remark 4.4** Note that the boundedness property of Theorem 4.2 (i) does not require explicitly the condition of Assumption 4.1.7 that $\mathbf{S} := \{T(s): s \in S\}$ is asymptotically nonexpansive. On the other hand, neither Theorem 4.2 nor Theorem 4.3 require Assumption 4.1.3. □

**Definition 4.5**, [8]. Let the sequence of means $\{\mu_k\}$ be in $Z \subset \ell_\infty(S)$ and let $\mathbf{S} := \{T(s): s \in S\}$ be a representation of a left reversible semigroup $S$. Then $Z$ is $\mathbf{S}$-stable if the functions $s \mapsto \langle T(s)x, x^* \rangle$ and $s \mapsto \|T(s)x - y\|$ on $S$ are also in $Z \subset \ell^\infty(S)$; $\forall x, y \in C$, $\forall x^* \in X^*$. □

**Definition 4.6**, [8, 11]. Let B and D be convex subsets of the Banach space X, with $\emptyset \ne D \subset B$ under proper inclusion, and let $P: B \to D$ be a retraction of B onto D. Then P is said to be sunny if $P(Px + t(x - Px)) = Px$; $\forall x \in B$, $\forall t \in \mathbf{R}_{0+}$ provided that $Px + t(x - Px) \in B$. □

**Definition 4.7**. D is said t be a sunny nonexpansive retract of B if there exists a sunny non expansive retraction P of B onto D. □

It is known that if C is weakly compact, $\mu$ is a mean on Z ( see Definition 3.5) and $s \mapsto \langle T(s)x, x^* \rangle$ is in Z for each $x^* \in X^*$ then there is a unique $x_0 \in X$ such that $\mu_s \langle T(s)x, x^* \rangle = \langle x_0, x^* \rangle$ for each $x^* \in X^*$. Also, if X is smooth, i.e. the duality mapping J of X is single-valued then a retraction P of B onto D is sunny and nonexpansive if and only if $\langle x - Px, J(z - Px) \rangle \le 0$; $\forall x \in B, \forall z \in D$ [6], [11].

**Remarks 4.8**. Note that Theorem 4.3 proves the convergence to a fixed point in **S**, with $F(\mathbf{S})$ being constructively proven to be nonempty by first building a sufficiently large convex compact $C_0$ so that the solution of the iterative scheme (4.1) is always bounded on $C_0$. Note also that Theorems 4.2 and 4.3 need not the assumption of $Z \subset \ell^\infty(S)$ being a left invariant $S$-stable subspace of containing "1" and to be a left-invariant mean on $Z$, although it is assumed to be strongly left regular so that so that it fulfils $\lim_{k \to \infty} \|\mu_{k+1} - \mu_k\| = 0$; $\forall \mu_k \in Z^*$ (Assumption 4.1.7)-see Definitions 3.6. However, the convergence to a unique fixed point in the set $F(S)$ is not proven under those less stringent assumptions. Note also that Assumption 4.1.8 required by Theorem 4.2 and also by Theorem 4.3 as a result is one of the two properties associated with the $S$-stability of $Z$. □

The results of Theorem 4.2 and Theorem 4.3 with further considerations by using Definitions 4.5-4.6 allow to obtain the convergence to a unique fixed point under more stringent conditions for the semigroup of self-mappings $T(\mu_k): C \to C$, $\mu_k \in Z^*$ as follows:

**Theorem 4.9**. If Assumptions 4.1 hold and, furthermore, $Z$ is a left-invariant $S$-stable subspace of $\ell^\infty(S)$ then the sequence $\{x_k\}$, generated by (4.1), converges strongly to a unique $x^* \in F(S)$; $\forall x_0 \in C$, $\forall \omega_k \in W$, $\forall r \in \mathbf{R}$ which is the unique solution of the variational inequality $\langle (f-I)x^*, J(y-x^*) \rangle \le 0$, $\forall y \in F(S)$. Equivalently, $x^* = P f x^*$ where $P$ is the unique sunny nonexpansive retraction of $C$ onto $F(S)$. □

The proof follows under similar tools as those used in [6] since $F(S)$ is a nonempy sunny nonexpansive retract of C which is unique since $T(\mu_k)$ is nonexpansive for all $\mu_k \in Z^*$.

Let $\{\bar{x}_k\}$ be the sequence solution generated by the particular iterative scheme resulting from (4.1) for any initial conditions $\bar{x}_0 = x_0 \in C$ when all the functions $\varphi_j$ and r are zeroed. It is obvious by the calculation of the recursive solution of (4.1) from (4.12) that the error from both solutions satisfies

$$\bar{x}_k = x_k - \sum_{\ell=0}^{k-1}\left(\prod_{j=\ell+1}^{k-1}[\beta_j]\right)\left(\sum_{j=1}^{s_\ell} v_{j\ell}\varphi_j(d(\omega_\ell, \omega_{\ell-p})) + \delta_\ell r\right) \; ; \; \forall k \in \mathbf{Z}_+$$

Since the convergence of the solution to fixed points of Theorems 4.2, 4.3 and 4.9 follows also for the sequence $\{\bar{x}_k\}$ it follows that a unique fixed point exists satisfying

$$\bar{x}^* = x^* - \sum_{\ell=0}^{\infty}\left(\prod_{j=\ell+1}^{\infty}[\beta_j]\right)\left(\sum_{j=1}^{s_\ell} v_{j\ell}\varphi_j(d(\omega_\ell, \omega_{\ell-p})) + \delta_\ell r\right)$$

where $\bar{x}^* \in F(S)$ is unique since $x^* \in F(S)$ is also unique from Theorem 4.9. Assume that $\beta_i \in (0, \beta)$ with $\beta < 1$. Then,

$$|\bar{x}^* - x^*| \le \lim_{k \to \infty}\sum_{\ell=0}^{k}\left(\prod_{j=\ell+1}^{k}[\beta_j]\right)\left(\sum_{j=1}^{s_\ell} v_{j\ell}\varphi_j(d(\omega_\ell, \omega_{\ell-p})) + \delta_\ell r\right) \le \frac{1}{1-\beta}\limsup_{k \to \infty}\left|\sum_{j=1}^{s_\ell} v_{j\ell}\varphi_j(d(\omega_\ell, \omega_{\ell-p})) + \delta_\ell r\right|$$

If $\delta_k = 1$ and $\beta_k = \beta < 1$; $\forall k \in \mathbf{Z}_{0+}$ and the $\varphi_j$-functions are zero then both fixed points are related by the constraint $\bar{x}^* = x^* - \dfrac{r}{1-\beta}$. Thus, consider a representation $\bar{\mathbf{S}} := \{T(s): s \in \bar{S}\}$ of a left reversible semigroup $\bar{S}$ as Lipschitzian mappings on $\bar{C}$ (see Definition 3.2 and Definition 3.3), a nonempty compact subset of de smooth Banach space X with Lipschitz constants $\{\bar{k}(s): s \in S\}$ which is asymptotically nonexpansive. Consider the iteration scheme:

$$\bar{x}_{k+1} = \beta_k x_k + \alpha_k f(\bar{x}_k) + \gamma_k T(\mu_k)\bar{x}_k = \beta_k x_k + (1-\beta_k)\bar{z}_k \quad (4.15)$$

$$\bar{z}_k = \dfrac{1}{1-\beta_k}(\alpha_k f(\bar{x}_k) + \gamma_k T(\mu_k)\bar{x}_k) \quad (4.16)$$

with $\bar{x}_0 \in \bar{C}$, where:

- $\{\mu_k\}$ is a strongly left regular sequence of means on $Z \subset \ell^\infty(S)$, i.e. $\mu_k \in Z^*$ (the dual of Z). See Definition 3.5 and Definition 3.6.

- $\bar{S}$ is a left reversible semigroup represented as Lipschitzian mappings on $\bar{C}$ by $\bar{\mathbf{S}} := \{T(s): s \in \bar{S}\}$.

□

**Assumptions 4.10**: The iterative scheme (4.15) keeps the applicable parts of Assumptions 4.1.1-4.1.5, 4.1.8 for the non identically zero parameterizing sequences $\{\alpha_k\}, \{\gamma_k\}$ and $\{\beta_k\}$. Assumptions 4.1.6 and 4.1.7 are modified with the replacements $C \to \bar{C}$, $S \to \bar{S}$ and $\mathbf{S} \to \bar{\mathbf{S}}$. □

Theorems 4.2 and 4.9 result in the following result for the iterative scheme (4.15) for $T(\mu_k): \bar{C} \to \bar{C}$, $\mu_k \in Z^*$:

**Theorem 4.11**. The following properties hold under Assumptions 4.10:

**(i)** $\max\left(\sup\limits_{k \in \mathbf{Z}_{0+}} |\bar{x}_k|, \sup\limits_{k \in \mathbf{Z}_{0+}} |T(\mu_k)\bar{x}_k|\right) < \infty$; $\forall \bar{x}_0 \in \bar{C}$. Also, $\|\bar{x}_k\| < \infty$ and $\|T(\mu_k)\bar{x}_k\| < \infty$ for any norm defined on the smooth Banach space X and there exists a nonempty bounded compact convex set $\bar{C}_0 \subseteq \bar{C} \subset X$ such that the solution of (4.2) is permanent in $\bar{C}_0$, $\forall k \geq k_0$ and some sufficiently large finite $k_0 \in \mathbf{Z}_{0+}$ with $\max\limits_{k \geq k_0}(\|\bar{x}_k\|, \|T(\mu_k)\bar{x}_k\|) \leq d_{\bar{C}_0} := \mathrm{diam}\,\bar{C}_0$.

**(ii)** $\lim\limits_{k \to \infty} \|T(\mu_k)\bar{x}_k - \bar{x}_k\| = 0$ and $\bar{x}_k \to \bar{z}_k \to \dfrac{\gamma_k T(\mu_k)\bar{x}_k}{1-\beta_k} \to T(\mu_k)\bar{x}_k \to \bar{x}^* \in \bar{C}_0$ as $k \to \infty$.

**(iii)** $\infty > \left|\bar{x}^* - \bar{x}_0\right| = \left|\lim\limits_{k \to \infty} \sum\limits_{j=0}^{k}(\bar{x}_{j+1} - \bar{x}_j)\right| = \left|\sum\limits_{j=0}^{\infty}(\alpha_j f(\bar{x}_j) + (\beta_j - 1)\bar{x}_j + \gamma_j T(\mu_j)\bar{x}_j)\right|$

**(iv)** Assume that the nonempty convex subset $\bar{C}$ of the smooth Banach space X, which contains the sequence $\{\mu_k\}$ of means on Z, is such that each element $\mu_k \in \mathbf{R}_{0+}^m$; $\forall k \in \mathbf{Z}_{0+}$, for some $m \in \mathbf{Z}_+$ so that (4.1) is a positive viscosity iteration scheme (4.15). Then,

**(iv.1)** $\{\bar{x}_k\}$ is a nonnegative sequence (i.e. all its components are nonnegative $\forall k \geq 0$, $\forall \bar{x}_0 \in \bar{C}$), denoted as $\bar{x}_k \geq 0$; $\forall k \geq 0$.

**(iv.2)** Property (i) holds for $\bar{C}_0 \subseteq \bar{C}$ and Property (ii) also holds for a limiting point $\bar{x}^* \in \bar{C}_0$.

**(iv.3)** Property (iii) becomes:

$$\infty > \left|\bar{x}^* - \bar{x}_0\right| = \left|\sum_{j=0}^{\infty} \left(\alpha_j f(\bar{x}_j) + \gamma_j T(\mu_j)\bar{x}_j\right) - \sum_{j=0}^{\infty} \left((1-\beta_j)\bar{x}_j\right)\right|$$

what implies that either

$$\sum_{j=0}^{\infty} \left(\alpha_j f(\bar{x}_j) + \gamma_j T(\mu_j)\bar{x}_j\right) < \infty \text{ and } \sum_{j=0}^{\infty} \left((1-\beta_j)x_j\right) < \infty$$

or

$$\limsup_{k \to \infty} \sum_{j=0}^{k} \left(\alpha_j f(\bar{x}_j) + \gamma_j T(\mu_j)\bar{x}_j\right) = \infty, \text{ and } \limsup_{k \to \infty} \sum_{j=0}^{\infty} \left((1-\beta_j)\bar{x}_j\right) = \infty$$

**(v)** If, furthermore, $Z$ is a left-invariant $\bar{S}$-stable subspace of $\ell^{\infty}(\bar{S})$ then the sequence $\{\bar{x}_k\}$, generated by (4.1), converges strongly to a unique $\bar{x}^* \in F(\bar{S})$; $\forall \bar{x}_0 \in \bar{C}$, $\forall r \in \mathbf{R}$ which is the unique solution of the variational inequality $\left\langle (f-I)\bar{x}^*, J(y-\bar{x}^*) \right\rangle \leq 0$, $\forall y \in F(\bar{S})$. Equivalently, $\bar{x}^* = \bar{P}f\bar{x}^*$ where $\bar{P}$ is the unique sunny nonexpansive retraction of $\bar{C}$ onto $F(\bar{S})$. Furthermore, the unique fixed points of the iterative schemes (4.1) and (4.15) are related by:

$$\bar{x}^* = x^* - \sum_{\ell=0}^{\infty} \left(\prod_{j=\ell+1}^{\infty} [\beta_j]\right) \left(\sum_{j=1}^{s_\ell} v_{j\ell} \varphi_j \left(d(\omega_\ell, \omega_{\ell-p})\right) + \delta_\ell r\right)$$

If, in addition, $\delta_k = 1$ and $\beta_k = \beta < 1$; $\forall k \in \mathbf{Z}_{0+}$ and the $\varphi_j$-functions are identically zero in the iterative scheme (4.1) then $\bar{x}^* = x^* - \dfrac{r}{1-\beta}$. □

**Remark 4.12**. Note that the results of Section 4 generalize those of Section 2 since the iterative process (4.1) possess simultaneously a nonlinear contraction and a nonexpansive mapping plus terms associated to driving terms combining both external driving forces plus the contribution of a nonlinear function evaluating distances over, in general, distinct metric spaces than that generating the solution of the iteration process. Therefore, the results about fixed points in Theorem 2.1 (vi)–(vii) are directly included in Theorem 4.2. □

Venter´s theorem can be used for the convergence to the equilibrium points of the solutions of the generalized iterative schemes (4.1) and (4.15), provided they are positive, as follows.

**Corollary 4.13**. Assume that:

1. $f, T(\mu_k): C \times \mathbf{Z}_{0+} \to \mathbf{R}_{0+}^m$ are both contractive mappings with $\varnothing \neq C \subset \mathbf{R}_{0+}^m$ being compact and convex, $\{\mu_k\}_{k \in \mathbf{Z}_{0+}} \in Z^*$, such that $Z$ is a left-invariant $\mathbf{S}$-stable subspace of $\ell^\infty(S)$ with S being a left-reversible semigroup.

2. $x_0, \bar{x}_0 \in C \subset \mathbf{R}_{0+}^m$, $r \in \mathbf{R}_{0+}$, with $C \supset \{0\}$ being compact and convex, $\alpha_k \in [0, \alpha]$, $\gamma_k \in [0, \gamma]$, $\delta_k \in [0, \delta]$ and $\beta_k \in [0,1)$; $\forall k \in \mathbf{Z}_{0+}$ for some real constants $\alpha, \gamma, \delta \in [0,1)$, and $\sum_{k=0}^{\infty} \delta_k < \infty$ if $r \neq 0$.

3. $\lim_{k \to \infty} \left( \sum_{j=0}^{k} (1-\beta_j) \right) = +\infty$ and $\exists \lim_{k \to \infty} \beta_k = 1$.

Then, the sets of fixed points of the positive iteration schemes (4.1) and (4.15) contain a common stable equilibrium point $0 \in \mathbf{R}_{0+}^m$ which is e unique solution to the variational equations of Theorem 4.9 and 4.11; i.e. $F(\mathbf{S}) \cap F(\bar{\mathbf{S}}) \supset \{0\}$ and that $x^* = \bar{x}^* = 0$.

**Outline of Proof**: The fact that the mappings $f, T(\mu_k): C \times \mathbf{Z}_{0+} \to \mathbf{R}_{0+}^m$ are both contractive, $\sum_{k=0}^{\infty} r \delta_k < \infty$ and $x_0, r \in \mathbf{R}_{0+}$ imply that the generated sequences $\{x_k\}, \{\bar{x}_k\}$ are both nonnegative and bounded for any $x_0, \bar{x}_0 \in C \subset \mathbf{R}_{0+}^m$ and they have unique zero limits from Theorem 2.1(v). □

The following result is obvious since if the representation S is nonexpansive, contractive or asymptotically contractive (Definition 3.3 and Definitions 3.4) then it is also asymptotically nonexpansive as a result:

**Corollary 4.14**. If the representation $\mathbf{S} := \{T(s) : s \in S\}$ is nonexpansive, contractive or asymptotically contractive then Theorem 4.2, Theorem 4.3 and Theorem 4.9 still hold under Assumptions 4.1, and Theorem 4.11 still holds under Assumptions 4.10. □

**ACKNOWLEDGMENTS**

The author is very grateful to the Spanish Ministry of Education by its partial support of this work through Grant DPI2006-00714. He is also grateful to the Basque Government by its support through GIC07143-IT-269-07.

**REFERENCES**

[1] T.A. Burton, *Stability by Fixed Point Theory for Functional Differential Equations*, Dover Publications Inc., Mineola, New York, 2006.

[2] P.N. Dutta and B.S. Choudhury, "A generalisation of contraction principle in metric spaces", *Fixed Point Theory and Applications*, Vol. 2008, Article ID 406638, 8 pages, doi:10.1155/2008/406368.

[3] C. Chen and C.Zhu, "Fixed point theorems for n times reasonable expansive mapping", *Fixed Point Theory and Applications*, Vol. 2008, Article ID 302617, 6 pages, doi:10.1155/2008/302617.

[4] L.G.Hu, "T-Strong convergence of a modified Halpern´s iteration for nonexpansive mappings", *Fixed Point Theory and Applications*, Vol. 2009, Article ID 418971, 8 pages, doi:10.1155/2009/418971.

[5] Y. Qing and B.E. Rhoades, "T-stability of Picard iteration in metric spaces", *Fixed Point Theory and Applications*, Vol. 2008, Article ID 418971, 4 pages, doi:10.1155/2008/418971.

[6] S. Saeidi, "Approximating common fixed points of Lipschitzian semigroup in smooth Banach spaces", *Fixed Point Theory and Applications*, Vol. 2009, Article ID 418971, 14 pages, doi:10.1155/2009/418971.

[7] J.M. Mendel, *Discrete Techniques of Parameter Estimation. The Equation Error Formulation*. <u>Control Theory: A series of Monographs and Textbooks</u> .Marcel Dekker Inc., New York, 1973.

[8] T. Suzuki, " Strong convergence of Kranoselskii and Mann´s type sequences for one-parameter nonexpansive semigroups without Bochner integrals", *Journal of Mathematical Analysis and Applications*, vol. 305, no. 1, pp. 227-239, 2005.

[9] B. Halpern, "Fixed points of nonexpansive maps", *Bull. Amer. Math. Soc.*, Vol. 73, pp. 957-961, 1967.

[10] P. L. Lions, "Approximation de points fixes de contractions ", *C.R. Acad. Sci. Paris, Ser. A-B*, Vol. 284, A1357-A1359, 1977.

[11] A. Aleyner and S. Reich, "An explicit construction of sunny nonexpansive retractions in Banach spaces ", *Fixed Point Theory and Applications*, Vol. 2005, No.3 pp. 295-305, 2005.

[12] M. De la Sen, "About robust stability of dynamic systems with time-delays through fixed point theory ", *Fixed Point Theory and Applications* (in press).

[13] A.T. Lau, H. Miyake and W. Takahashi, "Approximation of fixed points for amenable semigroups of nonexpansive mappings in Banach spaces", *Nonlinear Analysis*, Vol. 67, No. 5, pp. 111-1225, 2007.

[14] H. K. Xu, "Viscosity approximation methods for nonexpansive mappings", *Journal of Mathematical Analysis and Applications*, Vol. 298, No.1, pp. 279-291, 2004.

[15] A. Moudafi, "Viscosity approximation methods for fixed-points problems", *Journal of Mathematical Analysis and Applications*, Vol. 12, No. 3, pp. 46-55, 2000.